\documentclass[12pt]{article}
\usepackage{amsmath}
\usepackage{amssymb}
\usepackage{amscd}
\usepackage[T1]{fontenc}
\usepackage{xspace}

\newtheorem{theorem}{THEOREM}

\newtheorem{proposition}[theorem]{PROPOSITION}

\newtheorem{definition}{DEFINITION}

\newcommand{\G}{\Gamma}

\newcommand{\beq}{\begin{equation}\label}
\newcommand{\eeq}{\end{equation}}

\newcommand{\op}{\operatorname}
\newcommand{\Id}{\op{Id}}

\def\dim{{\mbox{dim}}}

\def\cala{{\cal A}}

\def\cale{{\cal E}}

\def\bbbone{\mbox{\rm 1\hspace {-.6em} l}}

\def\fracso{{\mathfrak {so}}}

\begin{document}

\baselineskip=0.7cm

\begin{center}
  \thispagestyle{empty}
{\large\bf INHOMOGENEOUS YANG-MILLS ALGEBRAS}
\end{center}
\vspace{0.75cm}

\begin{center}
  Roland BERGER
  \footnote{LaMUSE, Facult\'e des Sciences et Techniques,
  23 rue P. Michelon,
  F-42023 Saint-Etienne Cedex 2, France\\
Roland.Berger@univ-st-etienne.fr\\
  } and 
  Michel DUBOIS-VIOLETTE
\footnote{Laboratoire de Physique Th\'eorique, UMR 8627, Universit\'e Paris XI,
B\^atiment 210, F-91 405 Orsay Cedex, France\\
Michel.Dubois-Violette$@$th.u-psud.fr\\
},

\end{center} \vspace{1cm}

\begin{center} \today \end{center}

\vspace {1cm}

\begin{abstract}
We determine all  inhomogeneous Yang-Mills algebras and super Yang-Mills algebras which are Koszul. Following a recent proposal, a non-homogeneous algebra is said to be Koszul if the homogeneous part is Koszul  and if the PBW property holds. In this paper, the homogeneous parts are the Yang-Mills algebra and the super Yang-Mills algebra.
\end{abstract}

   \vspace{1,5cm}
\noindent LPT-ORSAY 05-72

\newpage

\section{Introduction}
After the introduction in \cite{ber:2001a} of the notion of Koszulity for $N$-homogeneous algebras, various concepts and constructions relative to quadratic algebras have been extended to $N$-homogeneous algebras \cite{ber-mdv-wam:2003}. An illustration of this, connected with theoretical physics, is provided by the study of the Yang-Mills algebra described in \cite{ac-mdv:2002b} and in \cite{ac-mdv:2004}. More recently, an extension to nonhomogeneous algebras of the notion of Koszulity has been proposed (and justified) in \cite{ber-gin:2005}. Besides the Koszulity for the homogeneous part, this notion involves a generalization of the Poincar\'e-Birkhoff-Witt property (PBW property) which has been also considered in \cite{flo-vat:2005}.\\

Let $\cala$ be the $N$-homogeneous algebra generated by $s+1$ elements $x^\lambda$ ($\lambda\in \{0,\dots,s\}$) with relations
\begin{equation}
W_{a\,\, \lambda_1\dots \lambda_N}x^{\lambda_1}\dots x^{\lambda_N}=0,\,\,\, a\in \{1,\dots,r\}
\end{equation}
then an {\sl inhomogeneous version of } $\cala$ is an algebra $\cala^J$, again generated by $s+1$ elements $x^\lambda$, but with relations
\begin{equation}
W_{a\,\, \lambda_1\dots \lambda_N}x^{\lambda_1}\dots x^{\lambda_N}=J_a,\,\,\, a\in \{1,\dots,r\}
\end{equation}
with 
\begin{equation}
J_a=\sum^{N-1}_{p=0}j_{a\,\, \lambda_1\dots \lambda_p} x^{\lambda_1}\dots x^{\lambda_p}
\end{equation}
The algebra $\cala$ is graded canonically (with $x^\mu$ of degree 1) while there is  only a filtration on $\cala^J$. Let $B$ be a basis of $\cala$ which consists of monomials in the generators and let $\tilde B$ be a corresponding family of monomials in the free algebra over the generators, then the classes in $\cala^J$ of these monomials generate the vector space $\cala^J$ and they form a basis if and only if the PBW property is satisfied.  In other words the PBW property means that the linear mapping of $\cala$ onto $\cala^J$ associated to $B$ and $\tilde B$ as above is an isomorphism. By introducing the associated graded algebra $gr(\cala^J)$ one obtains a canonical homomorphism of $\cala$ onto $gr(\cala^J)$ which is an isomorphism if and only if the PBW property is satisfied.\\

In the case of the Yang-Mills algebra, the generators are denoted by $\nabla_\lambda$ ($\lambda\in \{0,\dots,s\}$) and the homogeneous relations are $s+1$ cubic relations $W^\rho=0$ ($\rho\in \{0,\dots,s\}$) so for the inhomogeneous Yang-Mills algebras the relations (1.2) read $W^\rho=J^\rho$. One has identically $[\nabla_\rho,W^\rho]=0$ so the PBW property can be satisfied only if $[\nabla_\rho,J^\rho]=0$ follows from the relations (1.2). It turns out that this is also sufficient as will be shown in Section 3.\\

The plan of the paper is the following one. In the next section, Section 2, we review the results and the definition of \cite{ber-gin:2005} within an historical perspective. In Section 3, we apply the above results to determine all Koszul inhomogeneous Yang-Mills algebras in the sense of \cite{ber-gin:2005} and point out that the PBW property here is equivalent to the ``covariant conservation of the current" $[\nabla_\mu,J^\mu]=0$. In Section 4 we determine all inhomogeneous super Yang-Mills algebras. In both cases, the homogeneous parts were shown to be Koszul in \cite{ac-mdv:2002b} and \cite{ac-mdv:2004} respectively. The Yang-Mills algebra and its super version are not only Koszul algebras but are also Gorenstein. This latter property is a sort of Poincar\'e duality. At the Hochschild homological level, the sense in which it is a form of Poincar\'e duality has been made precise in \cite{ber-mar:2003}.\\

Throughout this paper $\mathbb K$ denotes a field and if $E$ and $F$ are $\mathbb K$-vector spaces, their tensor product over $\mathbb K$ is simply denoted by $E\otimes F$ while the tensor algebra, the symmetric algebra and the exterior algebra of $E$ are denoted respectively by $T(E)$, $S(E)$ and $\wedge(E)$. We use also everywhere Einstein's convention of summation of repeated up-down indices.

\section{Non-homogeneous Koszul algebras}
\setcounter{equation}{0}

Let $(V, \psi)$ be a Lie algebra, where $\psi : V\times V \rightarrow V$ is a Lie bracket over the $\mathbb K$ vector space $V$. The construction of the enveloping algebra $U_{\psi}$ of the Lie algebra $(V, \psi)$, and the PBW theorem are well-known. Denote by $(P)$ the two-sided ideal generated by $P$ for any subset $P$ of $T(V)$. Then 
$$U_{\psi} = T(V) / (\{x\otimes y-y \otimes x-\psi (x,y)\, ; \ x,\, y \in V\}).$$
The algebra $U_{\psi}$ is naturally filtered. Denote by $gr(U_{\psi})$ the associated graded algebra. The homogeneous quadratic part of the relations provides the symmetric algebra $S(V)=U_0$ which is also a graded algebra. The PBW theorem asserts that the natural graded algebra morphism $p: S(V) \rightarrow gr(U_{\psi})$ is an isomorphism. The PBW theorem admits a generalization to any inhomogeneous quadratic algebras \cite{bra-gai:1996} (see also \cite{pol-pos:1994}) which we now describe.

Consider the natural filtration $F^n= \bigoplus _{1\leq i \leq n} V^{\otimes i}$, $n=0,1,\ldots$, of $T(V)$. Let $P$ be a subspace of $F^2$. Introduce the filtered algebra $U=T(V) / (P)$. Set $R=\pi(P)\subseteq V\otimes V$, where $\pi$ is the natural projection of $F^2$ onto $V\otimes V$. The homogeneous quadratic part of the relations of $U$ provides the graded algebra $A=T(V) / (R)$. Let us say that $U$ satisfies the PBW property if the natural graded algebra morphism $p:A\rightarrow gr(U)$ is an isomorphism. 
\begin{theorem} \label{2PBW}
\emph{(PBW theorem in the quadratic case)} With the above notations and assumptions, suppose that the homogeneous quadratic algebra $A$ is Koszul.
Then $U$ satisfies the PBW property if and only if the conditions 
\begin{equation} \label{I}
P\cap F^{1} =0,
\end{equation}
\begin{equation} \label{J}
(PV+VP) \cap F^{2} \subseteq P
\end{equation}
hold.
\end{theorem}

If $U=U_{\psi}$ as at the beginning but $\psi : V\times V \rightarrow V$ is just an alternate bilinear map, then (\ref{I}) is clear and (\ref{J}) is equivalent to the Jacobi identity for $\psi$. So (\ref{J}) can be viewed as a generalization of the Jacobi identity. The proof of Theorem 1 given in \cite{bra-gai:1996} lies on Deformation Theory, while the proof in \cite{pol-pos:1994} is more elementary.

The \emph{symplectic reflection algebras} introduced in \cite{eti-gin:2002} satisfy a PBW theorem in a new context. This context is given by a finite group $\G$ acting on a vector space $V$. The symplectic reflection algebras are multi-parameter deformations of a certain basic algebra $S_{\bf{k}} (\cale)$ associated to $\G$ and $V$, and various specializations of the parameters (or of the basic algebra) provide several families of algebras which are of great interest in Representation Theory and Algebraic Geometry (desingularization). In this context, the ground field $\mathbb{K}$ is replaced by a certain semi-simple non-commutative ring $\bf{k}$, and we are concerned on deformations of the (non-commutative!) symmetric algebra $S_{\bf{k}} (\cale)$ of some $\bf{k}$-$\bf{k}$-bimodule $\cale$. Let us enter into more details. 

 Let $V$ be a finite dimensional complex vector space which is endowed with a symplectic 2-form $\omega$ (here $\mathbb{K}=\mathbb{C}$). Let $\G$ be a finite subgroup of Sp$(V)$. Denote by $\mathbf{k} =\mathbb{C}\G$ the group algebra of $\G$. Consider
the $\bf{k}$-$\bf{k}$-bimodule $\cale=V\otimes \bf{k}$, with left $\Gamma$-action given
by $g : v\otimes a \mapsto g(v)\otimes (ga)$, and right
$\Gamma$-action given by $v\otimes a \mapsto v\otimes (ag)$,
where $ga$ and $ag$ stand for the product in the group algebra. Then the $\bf{k}$-$\bf{k}$-algebras $T_{\bf{k}}(\cale)$ and $S_{\bf{k}}(\cale)$ can be considered. 

For any $g \in \G$, introduce the subspaces $M_g:=\op{Im}(\Id- g)$ and $L_g:=\op{Ker}(\Id- g)$. One has $V=M_g \oplus L_g$ and 
\begin{equation} \label{deco}
\wedge ^{2}(V)=(\wedge ^2(M_g)) \oplus (M_g\otimes L_g) \oplus (\wedge ^2(L_g)).
\end{equation}

The integer $a(g):=\dim M_g$ is even, and $g$ is called a \emph{symplectic reflection} if this dimension is 2. Let $\psi_g : \wedge ^{2}(V) \rightarrow \mathbb{C}$ be the $\mathbb{C}$-linear map which coincides with $\omega$ on $\wedge ^{a(g)}(M_g)\otimes \wedge ^{2-a(g)}(L_g)$ and vanishes on the other components of (\ref{deco}). Clearly, $\psi_g=0$ if $g$ is neither $Id$ nor a symplectic reflection. 

One defines a skew-symmetric $\G$-equivariant $\mathbb{C}$-bilinear map $\psi : V\times V \rightarrow \bf{k}$ by setting $\psi=\sum_{g\in\G}\psi_g\cdot g$. Then the \emph{symplectic reflection algebra}~\cite{eti-gin:2002} is the $\bf{k}$-algebra H$_{\psi}$ defined by
$$H_{\psi} = T_{\bf{k}}(\cale) / (\{x\otimes y-y \otimes x-\psi (x,y)\, ; \ x,\, y \in V\}).$$
For any map $m:\G \rightarrow \mathbb{C}$ which is constant on any conjugacy class, $m\cdot \psi$ is $\G$-equivariant, and H$_{m\cdot \psi}$ is also called a symplectic reflection algebra. The $m$'s are the parameters of the deformation. Thus the dimension of the parameter space (considered as a projective space) is the number of conjugacy classes of symplectic reflections in $\G$. Using an extension of Theorem \ref{2PBW} to semisimple ground rings, the following PBW theorem for symplectic reflection algebras was proved in \cite{eti-gin:2002}.
\begin{theorem} \label{EGPBW}
The natural graded algebra morphism 
\[
H_{0\cdot \psi}= S_{\bf{k}}(\cale) \rightarrow gr(H_{m\cdot \psi})
\]
is an isomorphism.
\end{theorem}

In \cite{ber-gin:2005}, this theorem was generalized to a new class of H$_{m\cdot \psi}$ which can be described as follows. Let $N$ be an integer with $2\leq N \leq \dim V$, $\G$ be a finite subgroup of GL$(V)$, and $\phi : \wedge^N (V)\rightarrow \mathbb C$ be a $\G$-invariant linear map (playing the role of $\omega$). There are an analogous decomposition (\ref{deco}) for $\wedge^{N}(V)$, and analogous definitions for the $\psi_g$'s, $\psi$, $m\cdot \psi$, and H$_{m\cdot \psi}$. In the relations of H$_{m\cdot \psi}$, the 2-tensor $x\otimes y-y \otimes x$ is replaced by any totally skew-symmetric tensor of $N$ variables, while $m\cdot \psi$ is applying on these variables. A first benefit of this situation (in particular, in the symplectic case) is that the parameter space may contain more conjugation classes, up to contain any conjugation class. 

The algebras of the new class are called \emph{higher symplectic reflection algebras}.
For the new class, it was known \cite{ber:2001a} that the undeformed algebra (i.e., corresponding to $m=0$, and whose relations are the antisymmetrizers of degree $N$) is Koszul when $\G$ is trivial, i.e., when $\bf{k}=\mathbb{C}$ (next the change of rings from $\mathbb{C}$ to any $\mathbb{C}\G$ is easy). Theorem 2 for the new class is proved in \cite{ber-gin:2005} by proving an $N$-version of the quadratic PBW theorem and using an argument based on a standard Koszul complex. The $N$-version of the quadratic PBW theorem is the following.
\begin{theorem} \label{NPBW} 
Assume that $\bf{k}$ is a von Neumann regular ring, $V$ is a
$\bf{k}$-$\bf{k}$-bimodule, $N$ is an integer $\geq 2$, and $P$ is a sub-$\bf{k}$-$\bf{k}$-bimodule of $F^N$, where
$F^n=\oplus_{0\leq i \leq n} V^{\otimes i}$ for any $n \geq 0$ (tensor product over $\mathbf{k}$). Set $U=T_{\mathbf{k}}(V)/(P)$ and
$A=T_{\mathbf{k}}(V)/(R)$, where $R=\pi(P)$ and $\pi$ is the projection of $F^N$ onto $V^{\otimes N}$ modulo
$F^{N-1}$. 
Assume that the graded left $\bf{k}$-module Tor$_3^A(\bf{k},\bf{k})$ is concentrated in degree $N+1$ (this property holds if $A$ is Koszul).
Then $U$ satisfies the PBW property if and only if the conditions
\begin{equation} \label{NI}
P\cap F^{N-1} =0,
\end{equation}
\begin{equation} \label{NJ}
(PV+VP) \cap F^{N} \subseteq P
\end{equation}
hold.
\end{theorem}

Let us recall the definition of a von Neumann regular ring \cite{goo:1979}, \cite{wei:1994}. A ring $\bf{k}$ is said to be \emph{von Neumann regular} if for any $x\in\bf{k}$, there exists $y\in\bf{k}$ such that $xyx=x$. It is a basic fact that $\bf{k}$ is von Neumann regular if and only if any left (right) $\bf{k}$-module is \emph{flat}, and it is this fact which is used in the proof of Theorem \ref{NPBW}. Since any projective module is flat, any semi-simple algebra is von Neumann regular. Thus this theorem applies to higher symplectic reflection algebras ($\mathbb{C}\G$ is semi-simple). The proof of Theorem \ref{NPBW} given in \cite{ber-gin:2005} follows along the lines of the proof of Theorem \ref{2PBW} which is given in \cite{pol-pos:1994} and which is more adapted (via flatness) to an extension to von Neumann regular ground rings. 

In Theorem \ref{NPBW}, the condition which generalizes the Jacobi identity is (\ref{NJ}). This condition is expressed in a more explicit manner as follows.
\begin{proposition} \label{car}
Assume that (\ref{NI}) holds. Let $\varphi :R \rightarrow F^{N-1}$ be the $\bf{k}$-$\bf{k}$-linear map such
that $P=\{x-\varphi(x);\ x \in R \}$. Decompose $\varphi = \sum_{j=0}^{N-1}\varphi_j$, $\varphi_j : R\rightarrow
V^{\otimes j}$. Set $\mathcal{W}_{N+1}= (R\otimes V)\cap (V\otimes R)$, and denote by $I$ the identity map of $V$. Then (\ref{NJ}) is equivalent to the set of the following relations
\begin{equation} \label{J'1}
\left(\varphi_{N-1}\otimes I-I\otimes \varphi_{N-1}\right)(\mathcal{W}_{N+1})\subseteq R,
\end{equation}
\begin{equation} \label{J'2}
\begin{array}{l}
\left(\varphi_j(\varphi_{N-1}\otimes I-I\otimes \varphi_{N-1})+\varphi_{j-1}\otimes I-I\otimes \varphi_{j-1}
\right)(\mathcal{W}_{N+1})=0,\\
\\
 1\leq j \leq N-1,
\end{array}
\end{equation}
\begin{equation} \label{J'3}
\varphi_0(\varphi_{N-1}\otimes I - I\otimes \varphi_{N-1})(\mathcal{W}_{N+1})=0.
\end{equation}
\end{proposition}

An $N$-PBW theorem has been also obtained in \cite{flo-vat:2005}. It is Theorem \ref{NPBW} in which $\bf{k}$ is a field, $V$ is a finite dimensional vector space over $\bf{k}$, and the assumption on $A$ is the full one saying that $A$ is Koszul. The proof uses Deformation Theory and follows \cite{bra-gai:1996}. In \cite{flo-vat:2005}, the $N$-PBW theorem has been applied to a classification result determining the PBW deformations of the $N$-homogeneous algebras whose relations are the antisymmetrizers of degree $N$. The intersection of this classification result with higher symplectic reflection algebras is very small since it corresponds to a trivial group $\G$.

In the next section, we shall apply the $N$-PBW theorem to determine the PBW-deformations of the Yang-Mills algebras. It will be a situation in which the $N$-PBW theorem of \cite{flo-vat:2005} suffices. However we believe it is important to show that a general setting of the $N$-PBW theorem as developed in \cite{ber-gin:2005} allows to include significant ``relative'' situations where the ground field is enlarged to certain non-commutative rings. 

Note also that the PBW-deformations of the Yang-Mills algebras \emph{having two generators} (i.e., when $s=1$) will coincide with some special cases already obtained in \cite{flo-vat:2005} by applying the $N$-PBW theorem to the classification of the PBW-deformations of cubic AS-regular algebras of global dimension 3. For a higher number of generators, the Yang-Mills algebras are not AS-regular since they have an exponential growth (as graded algebras). 

Let us end this section by a terminology point which deserves to be mentioned. Besides the terminology ``PBW-deformations'' which is natural when the starting point is $A$ and the aim is $U$, it may be convenient to have an intrinsic terminology for the filtered algebras $U$. It was proposed in \cite{ber-gin:2005} by saying that $U$ is \emph{Koszul}. Let us recall precisely this definition from \cite{ber-gin:2005}, in which a historical argument for this terminology is also given (going back to the first Lie theory use by J. L. Koszul of his complex). Note that this terminology agrees with the graded one.  
\begin{definition} \label{koszul}
Assume that $\bf{k}$ is a von Neumann regular ring, $V$ is a
$\bf{k}$-$\bf{k}$-bimodule, $N$ is an integer $\geq 2$, and $P$ is a sub-$\bf{k}$-$\bf{k}$-bimodule of $F^N$, where
$F^n=\oplus_{0\leq i \leq n} V^{\otimes i}$ for any $n \geq 0$. Set $U=T_{\mathbf{k}}(V)/(P)$ and
$A=T_{\mathbf{k}}(V)/(R)$, where $R=\pi(P)$ and $\pi$ is the projection of $F^N$ onto $V^{\otimes N}$ modulo $F^{N-1}$. Then $U$ is said to be Koszul if the graded algebra (with
$N$-homogeneous relations) $A$ is Koszul and if the PBW property
holds for $U$.
\end{definition}
\section{Regular currents for the Yang-Mills algebra}
\setcounter{equation}{0}

Let $g_{\lambda \mu}$ be the components of a nondegenerate symmetric bilinear form on $\mathbb K^{s+1}$ in the canonical basis and let us denote by $g^{\lambda\mu}$ be the components of the nondegenerate symmetric bilinear form on the dual vector space of $\mathbb K^{s+1}$ such that $g_{\lambda\mu}g^{\mu\nu}=\delta^\nu_\lambda$ (i.e. given by the inverse matrix). Recall that the {\sl Yang-Mills algebra} \cite{ac-mdv:2002b} is the cubic algebra $\cala$ generated by $s+1$ elements $\nabla_\lambda$ ($\lambda\in \{0,\dots,s\}$) with relations
\begin{equation}
g^{\lambda\mu}[\nabla_\lambda,[\nabla_\mu,\nabla_\nu]]=0,\,\,\, \nu\in\{0,\dots,s\}
\end{equation}
that is $\cala=A(E,R)$ with $E=\oplus_\lambda\mathbb K\nabla_\lambda$ and $R\subset E^{\otimes^3}$ given by
\begin{equation}
R=\sum_\rho \mathbb K W^\rho=\sum_\rho \mathbb K g^{\lambda\mu} g^{\nu\rho}[\nabla_\lambda[\nabla_\mu,\nabla_\nu]_\otimes]_\otimes
\end{equation}
where the notations of \cite{ber-mdv-wam:2003} have been used. The following result was proved in \cite{ac-mdv:2002b}.
\begin{theorem}\label{YM}
The cubic Yang-Mills algebra $\cala$ is Koszul of global dimension 3 and is Gorenstein.
\end{theorem}
The basis $W^\rho$ ($\rho\in \{0,\dots,s\}$) of $R$ reads
\begin{equation}
W^\rho=W^{\rho \lambda \mu\nu}\nabla_\lambda \otimes \nabla_\mu \otimes \nabla_\nu
\end{equation}
with
\begin{equation}
W^{\rho\lambda\mu\nu}=g^{\rho\lambda}g^{\mu\nu}+g^{\rho\nu}g^{\lambda\mu}-2g^{\rho\mu}g^{\lambda\nu}
\end{equation}
which satisfies
\begin{equation}
W^{\lambda\mu\nu\rho}=W^{\rho\lambda\mu\nu}
\end{equation}
that is one has
\begin{equation}
w:=W^\rho\otimes \nabla_\rho=\nabla_\rho\otimes W^\rho
\end{equation}
and $w$ spans the space ${\cal W}_{N+1}={\cal W}_4$ of last section which is  an 1-dimensional subspace of $E^{\otimes^4}$. Using the concepts introduced in \cite{mdv:2005} one has $\cala=\cala(w,3)$ and $Q_w=\bbbone$.\\
Besides invariance by cyclic permutation (3.5), $w$ has another important property namely
\begin{equation}
W^{\rho\lambda\mu\nu}+W^{\rho\nu\lambda\mu}+W^{\rho\mu\nu\lambda}=0
\end{equation}
which is useful for the computations.\\

Our aim now is to apply Proposition \ref{car} to construct all linear mappings $\varphi:R\rightarrow \mathbb K \bbbone \oplus E\oplus E^{\otimes^2}$ such that the algebra $\cala_\varphi=T(E)/ (\{x-\varphi(x)\vert x\in R\})$ is Koszul in the sense of Definition 1, i.e. satisfies the PBW property.\\
Define the {\sl current} $J$ with components $J^\rho$ ($\rho\in \{0,\dots,s\}$) by
\begin{equation}
J^\rho=\varphi(W^\rho)=j^{\mu\nu\rho} \nabla_\mu\otimes \nabla_\nu+j^{\lambda\rho} \nabla_\lambda + j^\rho\bbbone
\end{equation}
where the $j's$ are in $\mathbb K$. The current $J$ will be said to be {\sl regular} whenever $\cala_\varphi$ satisfies the PBW property (i.e. whenever $\cala_\varphi$ is Koszul).
 
 \begin{theorem}\label{IYM}
 The PBW property is satisfied, i.e. $J$ is regular, if and only if one has
 \[
 \begin{array}{lll}
 j^{\alpha\beta\gamma}& =&(g^{\alpha\rho}g^{\beta\gamma}-g^{\alpha\gamma}g^{\beta\rho})b_\rho+\omega^{\alpha\beta\gamma}+s^{\alpha\beta\gamma}\\
j^{\alpha\beta} & = & -\frac{1}{2}\omega^{\alpha\beta\rho}b_\rho+s^{\alpha\beta}\\
j^\alpha&=&s^\alpha
\end{array}
\]
where  $\omega^{\alpha\beta\gamma}$ is completely antisymmetric ($\omega\in \wedge^3\mathbb K^{s+1}$), $s^{\alpha\beta\gamma}$ is completely symmetric with $s^{\alpha\beta\rho}b_\rho=0$, $s^{\alpha\beta}$ is symmetric with $s^{\alpha\rho}b_\rho=0$ and $s^\alpha$ is such that $s^\rho b_\rho=0$.
 \end{theorem}
\noindent \underbar{Proof}. Equation (2.6) reads here
\[
j^{\alpha\beta\gamma}(\nabla_\alpha\otimes \nabla_\beta\otimes \nabla_\gamma-\nabla_\gamma\otimes \nabla_\alpha\otimes \nabla_\beta)=W^\rho b_\rho\,\, \mbox{ for some}\,\, b_\rho\in \mathbb K\, \,  (\rho\in \{0,\dots,s\})
\]
 which is equivalent to $j^{\alpha\beta\gamma}=(g^{\alpha\rho}g^{\beta\gamma}-g^{\beta\rho}g^{\alpha\gamma})b_\rho + c^{\alpha\beta\gamma}$ with $c^{\alpha\beta\gamma}=c^{\gamma\alpha\beta}$ ($\forall \alpha;\beta, \gamma \in \{0,\dots,s\}$) which means $c^{\alpha\beta\gamma}=\omega^{\alpha\beta\gamma}+s^{\alpha\beta\gamma}$ with $\omega^{\alpha\beta\gamma}$ completely antisymmetric and $s^{\alpha\beta\gamma}$ completely symmetric. Indeed $(g^{\alpha\rho}g^{\beta\gamma}-g^{\beta\rho}g^{\alpha\gamma})b_\rho$ is a solution of the above equation while $c^{\alpha\beta\gamma}$ invariant by cyclic permutation is the general solution of the corresponding equation for $b_\rho=0$.\\
 Equation (2.7) for $j=2$ reads here $(j^{\alpha\beta\rho}b_\rho+j^{\alpha\beta}-j^{\beta\alpha})\nabla_\alpha \otimes \nabla_\beta=0$ which is equivalent to $s^{\alpha\beta\rho}b_\rho=0$ and $j^{\alpha\beta}=-\frac{1}{2}\omega^{\alpha\beta\rho}b_\rho+s^{\alpha\beta}$ with $s^{\alpha\beta}=s^{\beta\alpha}$.\\
 Equation (2.7) for $j=1$ reads then (since $j^\rho\nabla_\rho-\nabla_\rho j^\rho=0$) $j^{\alpha\rho}b_\rho=0$ which is equivalent to $s^{\alpha\rho}b_\rho=0$. \\
 Finally, Equation (2.8) reads $j^\rho b_\rho=0$ which implies the result.$\square$\\
 
 The inhomogeneous Yang-Mills algebra $\cala_\varphi$ is the algebra generated by elements $\nabla_\lambda$ ($\lambda \in \{0,\dots,s\}$) with relations
 \begin{equation}
 [\nabla_\lambda, F^{\lambda\mu}]=J^\mu,\,\,\, \mu\in\{0,\dots,s\}
 \end{equation}
 with $F^{\mu\nu}=[\nabla^\mu,\nabla^\nu]=g^{\mu\alpha}g^{\nu\beta}[\nabla_\alpha,\nabla_\beta]=g^{\mu\alpha} g^{\nu\beta}F_{\alpha\beta}$ and where, with an abuse of notation, we have again denoted by $J^\mu$ the element of $\cala_\varphi$ obtained by multiplication in $\cala_\varphi$ from $J^\mu\in T(E)$ (i.e. by replacing the tensor products by the products in $\cala_\varphi$). On the other hand, the left hand side of (3.9) satisfies identically
 \begin{equation}
 [\nabla_\mu[\nabla_\lambda,F^{\lambda\mu}]]=0
 \end{equation}
 in view of the associativity of the product of the $\nabla_\alpha$. It follows that the PBW property can be satisfied only if one has
 \begin{equation}
 [\nabla_\mu,J^\mu]=0
 \end{equation}
 i.e. ``covariant conservation" of the current $J$. It turns out that this is sufficient, that is Equation (3.11) implies the conditions of the last theorem for $J$. Thus (3.11) is equivalent here to the PBW property and plays the role of the generalization of the Jacobi identity.\\
 To make contact with physics one should require that besides (3.11), $J$ only depends on the $\nabla_\lambda$ through the commutators $F_{\lambda\mu}=[\nabla_\lambda,\nabla_\mu]$ which implies
 \begin{equation}
 J^\mu=b_\lambda F^{\lambda\mu}  + \omega^{\lambda\rho\mu} F_{\lambda\rho}+s^\mu \bbbone
 \end{equation}
 with $\omega^{\alpha\beta\gamma}$ completely antisymmetric such that
 \begin{equation}
 b_\lambda \omega^{\lambda\mu\nu}=0
 \end{equation}
 and with $s^\alpha$ such that
 \begin{equation}
 b_\lambda s^\lambda=0
 \end{equation}
 in view Theorem (\ref{IYM}).\\
 Thinking of the electromagnetic case, $\nabla_\lambda=\partial_\lambda + A_\lambda$, etc., one sees that (3.12) completed with (3.13) and (3.14) contains a generalization of Ohm's law.
 
 \section{Inhomogeneous super Yang-Mills algebras}
 \setcounter{equation}{0}
 
 The Yang-Mills algebra is the universal enveloping algebra of a Lie algebra \cite{ac-mdv:2002b} which has a super variant. The {\sl super Yang-Mills algebra} defined in \cite{ac-mdv:2004} is the cubic algebra $\tilde\cala$ generated by $s+1$ elements $S_\lambda$ ($\lambda\in \{0,\dots,s\}$) with relations
 \begin{equation}
 g^{\lambda\mu}[S_\lambda,\{S_\mu,S_\nu\}]=0,\,\,\, \nu\in \{0,\dots,s\}
 \end{equation}
 where $g^{\alpha\beta},g_{\alpha\beta}$, etc. are as in the last section and where $\{A,B\}=AB+BA$. Relations (4.1) can be equivalently written as
 \begin{equation}
 [g^{\lambda\mu}S_\lambda S_\mu,\,\, S_\nu]=0,\,\,\, \nu\in\{0,\dots,s\}
 \end{equation}
 which means that the quadratic element $g^{\lambda\mu}S_\lambda S_\mu$ is a central element of $\tilde\cala$. Using again the notations of \cite{ber-mdv-wam:2003} one has $\tilde \cala=A(\tilde E,\tilde R)$ with $\tilde E=\oplus_\lambda \mathbb K S_\lambda$ and
 \begin{equation}
 \tilde R=\sum_\rho \mathbb K \tilde W^\rho = \sum_\rho \mathbb K(g^{\lambda\rho}g^{\mu\nu}-g^{\lambda\mu}g^{\nu\rho})S_\lambda\otimes S_\mu\otimes S_\nu
 \end{equation}
 The following result was proved in \cite{ac-mdv:2004}.
 \begin{theorem}
\label{SYM}
The cubic super Yang-Mills algebra $\tilde \cala$ is Koszul of global dimension 3 and is Gorenstein.
\end{theorem}
The basis $\tilde W^\rho\,\,\, (\rho\in\{,\dots,s\})$ of $\tilde R$ reads
\begin{equation}
\tilde W^\rho=\tilde W^{\rho\lambda\mu\nu} S_\lambda \otimes S_\mu \otimes S_\nu
\end{equation}
with 
\begin{equation}
\tilde W^{\rho\lambda\mu\nu}=g^{\rho\lambda}g^{\mu\nu}- g^{\rho\nu}g ^{\lambda\mu}
\end{equation}
which satisfies
\begin{equation}
\tilde W^{\lambda\mu\nu\rho}=-\tilde W^{\rho\lambda\mu\nu}
\end{equation}
so one has
\begin{equation}
\tilde w:=S_\rho \otimes \tilde W^\rho=-\tilde W^\rho\otimes S_\rho
\end{equation}
 $\tilde \cala=\cala(\tilde w,3)$ with $Q_{\tilde w}=-\bbbone$ (with the notations of \cite{mdv:2005}) and $\tilde w$ spans the space ${\cal W}_{N+1}={\cal W}_4$ of Section 2 which is again 1-dimensional.\\

We now add inhomogeneous terms
\[
\tilde J^\rho=\tilde \varphi(\tilde W^\rho)=\tilde \jmath^{\alpha\beta\rho}S_\alpha\otimes S_\beta+\tilde \jmath^{\alpha\rho}S_\alpha+ \tilde\jmath^\rho\bbbone\in \mathbb K\bbbone \oplus \tilde E\oplus \tilde E^{\otimes^2}
\]
and look for the conditions under which the PBW property holds for $\tilde\cala_{\tilde \varphi}=T(\tilde E)/(\{x-\tilde \varphi(x)\vert x\in \tilde R\})$ i.e. such that $\tilde \cala_{\tilde\varphi}$ is Koszul in the sense of \cite{ber-gin:2005}.

\begin{theorem}\label{ISYM}
The PBW property is satisfied if and only if one has
\[
\begin{array}{lll}
\tilde \jmath^{\alpha\beta\gamma} & = & (g^{\alpha\gamma}g^{\beta\rho}-g^{\beta\gamma}g^{\alpha\rho})b_\rho\\
\tilde\jmath^{\alpha\beta} & = & \omega^{\alpha\beta}\\
\tilde\jmath^\alpha & = & -\frac{1}{2} \omega^{\alpha\rho}b_\rho
\end{array}
\]
with $\omega^{\alpha\beta} =-\omega^{\beta\alpha}$.
\end{theorem}
\noindent \underbar{Proof}. Equation (2.6) reads here
$\tilde \jmath^{\alpha\beta\gamma}+\tilde\jmath^{\beta\gamma\alpha}=\tilde W^{\alpha\beta\gamma\rho} b_\rho=(g^{\alpha\beta}g^{\gamma\rho}-g^{\beta\gamma}g^{\alpha\rho})b_\rho$
 for some $b_\rho\in \mathbb K$  ($\rho\in \{0,\dots,s\}$) which is equivalent to 
 \begin{equation}
 \tilde\jmath^{\alpha\beta\gamma}=(g^{\alpha\gamma}g^{\beta\rho}-g^{\beta\gamma}g^{\alpha\rho})b_\rho
 \label{jtilde3}
 \end{equation}
 Equation (2.7) for $j=2$ reads here $\tilde\jmath^{\alpha\beta\gamma}b_\gamma+\tilde\jmath^{\alpha\beta}+\tilde \jmath^{\beta\alpha}=0$ which is equivalent to $\tilde\jmath^{\alpha\beta}=-\tilde\jmath^{\beta\alpha}=\omega^{\alpha\beta}$ since it follows from (\ref{jtilde3}) that one has $\tilde\jmath^{\alpha\beta\gamma}b_\gamma=0$ identically.\\
Equation (2.7) for $j=1$ reads here $\tilde\jmath^{\alpha\gamma}b_\gamma+2\tilde\jmath^\alpha=0$ which is equivalent to $\tilde\jmath^\alpha=-\frac{1}{2} \omega^{\alpha\rho}b_\rho$, while Equation (2.8) which reads $\tilde \jmath^\rho b_\rho=0$ is then identically satisfied. $\square$\\

With $\tilde J$ as in Theorem \ref{ISYM} above, the relations of the algebra read
\begin{equation}
[g^{\lambda\mu}S_\lambda S_\mu,S_\nu]=-[g^{\lambda\mu}b_\lambda S_\mu,S_\nu]+\omega^{\lambda\rho}g_{\rho\nu}(S_\lambda+\frac{1}{2}b_\lambda\bbbone)
\end{equation}
or, equivalently
\begin{equation}
[g^{\lambda\mu}\hat S_\lambda\hat S_\mu,\hat S_\nu]=\omega^{\tau\rho}g_{\rho\nu} \hat S_\tau
\end{equation}
with
\begin{equation}
\hat S_\lambda=S_\lambda+\frac{1}{2} b_\lambda \bbbone
\end{equation}
so $\mbox{ad}(g^{\lambda\mu}\hat S_\lambda \hat S_\mu$) induces the infinitesimal rotation $(\omega^{\tau\rho}g_{\rho\nu}) \in \fracso(s+1,g)$ in the generators $\hat S_\alpha$.

\end{document}